\documentclass[a4paper,10pt]{article}
\usepackage{times,url}
\usepackage{amsmath,amssymb,amsthm,amsfonts}

\newtheorem{theorem}{Theorem}[section]
\newtheorem{lemma}{Lemma}[section]
\newtheorem{corollary}{Corollary}[section]

\newtheorem{remark}{Remark}[section]

\begin{document}
\setcounter{page}{1} 

\begin{center}
{\LARGE \bf  On a Curious Biconditional Involving \\
Divisors of Odd Perfect Numbers}
\vspace{8mm}

{\large \bf Jose Arnaldo B. Dris}
\vspace{3mm}

PhD Mathematics Student, University of the Philippines-Diliman \\ 
e-mails: \url{jadris@feu.edu.ph}, \url{josearnaldobdris@gmail.com}
\vspace{2mm}

\end{center}
\vspace{10mm}

\noindent
{\bf Abstract:} We investigate the implications of a curious biconditional involving divisors of odd perfect numbers, if Dris conjecture that $q^k < n$ holds, where $q^k n^2$ is an odd perfect number with Euler prime $q$.  We then show that this biconditional holds unconditionally.  Lastly, we prove that the inequality $q<n$ holds unconditionally. 
\vspace{10mm}

\section{Introduction}
If $N$ is a positive integer, then we write $\sigma(N)$ for the sum of the divisors of $N$.  A number $N$ is \emph{perfect} if $\sigma(N)=2N$.  We denote the abundancy index $I$ of the positive integer $w$ as $I(w) = \sigma(w)/w$.  We also denote the deficiency $D$ of the positive integer $x$ as $D(x) = 2x - \sigma(x)$ \cite{OEIS-A033879}.

Euclid and Euler showed that an even perfect number $E$ must have the form
$$E=\left(2^p - 1\right){2^{p-1}}$$
where $2^p - 1$ is a \emph{Mersenne prime}.  On the other hand, Euler showed that an odd perfect number $O$ must have the form
$$O=q^k n^2$$
where $q$ is an \emph{Euler prime} (i.e., $q \equiv k \equiv 1 \pmod 4$ and $\gcd(q,n)=1$).

It is currently unknown whether there are any odd perfect numbers.  On the other hand, only $49$ even perfect numbers have been found, a couple of which were discovered by the Great Internet Mersenne Prime Search \cite{GIMPS}.  It is conjectured that there are infinitely many even perfect numbers, and that there are no odd perfect numbers.

Descartes, Frenicle and subsequently Sorli conjectured that $k=1$ \cite{Beasley}.  Sorli conjectured $k=1$ after testing large numbers with eight distinct prime factors for perfection \cite{Sorli}.

Dris conjectured in \cite{Dris1} and \cite{Dris2} that the divisors $q^k$ and $n$ are related by the inequality $q^k < n$.  Brown \cite{Brown} and Starni \cite{Starni} have recently uploaded preprints claiming a proof for the weaker inequality $q < n$.

Holdener presented some conditions equivalent to the existence of odd perfect numbers in \cite{Holdener}.  In this paper, we prove the following results:

\begin{lemma}
\label{lem:1}
If $N=q^k n^2$ is an odd perfect number with Euler prime $q$, then the sum
$$\frac{\sigma(q^k)}{n}+\frac{\sigma(n)}{q^k}$$
is bounded from above if and only if the sum
$$\frac{q^k}{n}+\frac{n}{q^k}$$
is bounded from above.
\end{lemma}

The following lemma is proved in the preprint \cite{DagalDris}.  (We will not need to use this result in the present paper.  Hence, we will not be proving this lemma here.)

\begin{lemma}
\label{lem:2}
If $N=q^k n^2$ is an odd perfect number with Euler prime $q$ and $3 \nmid N$, then $\sigma(q) \neq \sigma(n)$.
\end{lemma}

Using Lemma \ref{lem:1}, we are able to prove the following unconditional result.

\begin{theorem}
\label{thm:1}
If $N=q^k n^2$ is an odd perfect number with Euler prime $q$, then $\sigma(q^k) \neq \sigma(n)$.
\end{theorem}

\begin{lemma}
\label{lem:3}
If $N=q^k n^2$ is an odd perfect number with Euler prime $q$, then the inequality
$$\frac{\sigma(q^k)}{n}+\frac{\sigma(n)}{q^k}<I(q^k)+I(n)$$
holds if and only if the biconditional
$$q^k < n \iff  \sigma(n) < \sigma(q^k)$$
holds.
\end{lemma}

The following result is trivial.  The proof is easy, and is left for the interested reader.

\begin{lemma}
\label{lem:4}
If $N=q^k n^2$ is an odd perfect number with Euler prime $q$, then either $q^k < n$, $\sigma(q^k) < n$ or $\sigma(n) < q^k$ imply that the biconditional
$$q^k < n \iff \sigma(q^k) < \sigma(n) \iff \frac{\sigma(q^k)}{n} < \frac{\sigma(n)}{q^k}$$
holds.
\end{lemma}

The following corollary follows easily from Theorem \ref{thm:1} and Lemma \ref{lem:3}.

\begin{corollary}
\label{cor:1}
If $N=q^k n^2$ is an odd perfect number with Euler prime $q$, then the biconditional
$$q^k < n \iff \sigma(q^k) < \sigma(n) \iff \frac{\sigma(q^k)}{n} < \frac{\sigma(n)}{q^k}$$
holds.
\end{corollary}

All of the proofs given in this note are elementary.

\section{Preliminaries}
Let $N=q^k n^2$ be an odd perfect number with Euler prime $q$.

First, we show that the following equations hold. (The proof is taken from the paper \cite{Dris4}.)  This will serve as motivation for trying to prove the inequality $q^k < n$ or the stronger inequality $\sigma(q^k)<n$.

\begin{lemma}
\label{lem:5}
If $N=q^k n^2$ is an odd perfect number with Euler prime $q$, then
$$\gcd\left(n^2, \sigma(n^2)\right) = \frac{D(n^2)}{\sigma(q^{k-1})} = \frac{\sigma(N/q^k)}{q^k}.$$
\end{lemma}

\begin{proof}
Since $N=q^k n^2$ is an odd perfect number, we have
$$\sigma(q^k)\sigma(n^2) = \sigma(N) = 2N = 2{q^k}{n^2},$$
from which it follows that $q^k \mid \sigma(n^2)$ (because $\gcd\left(q^k,\sigma(q^k)\right)=1$). Hence,
$$\frac{\sigma(n^2)}{q^k} = \frac{\sigma(N/q^k)}{q^k}$$
is an integer.

First, we prove that
$$\frac{D(n^2)}{\sigma(q^{k-1})} = \frac{\sigma(N/q^k)}{q^k}.$$
We rewrite the equation
$$\sigma(q^k)\sigma(n^2) = 2{q^k}{n^2}$$
as
$$\left(q^k + \sigma(q^{k-1})\right)\sigma(n^2) = 2{q^k}{n^2}$$
$$\sigma(q^{k-1})\sigma(n^2) = {q^k}\left(2n^2 - \sigma(n^2)\right) = {q^k}\cdot{D(n^2)}$$
$$\frac{\sigma(n^2)}{q^k} = \frac{D(n^2)}{\sigma(q^{k-1})},$$
and we are done.

Next, we show that
$$\gcd\left(n^2, \sigma(n^2)\right) = \frac{D(n^2)}{\sigma(q^{k-1})}.$$
We already know that
$$\sigma(n^2) = {q^k}\cdot{\left(\frac{D(n^2)}{\sigma(q^{k-1})}\right)}.$$
Since $\sigma(q^k)\sigma(n^2) = 2{q^k}{n^2}$, we also obtain
$$\frac{2n^2}{\sigma(q^k)} = \frac{\sigma(n^2)}{q^k} = \frac{D(n^2)}{\sigma(q^{k-1})}.$$
This implies that
$$n^2 = \frac{\sigma(q^k)}{2}\cdot{\left(\frac{D(n^2)}{\sigma(q^{k-1})}\right)}.$$
It follows that
$$\gcd\left(n^2, \sigma(n^2)\right) = \frac{D(n^2)}{\sigma(q^{k-1})}$$
since
$$\gcd\left(q^k, \frac{\sigma(q^k)}{2}\right) = \gcd(q^k, \sigma(q^k)) = 1.$$
This concludes the proof.
\end{proof}

\begin{remark}
\label{rem:1}
Dris obtained the lower bound $3$ for $\sigma(N/q^k)/q^k$ in \cite{Dris1} and \cite{Dris2}.
\end{remark}

\begin{remark}
\label{rem:2}
Notice that
$$\frac{\sigma(n^2)}{q^k} = \frac{2n^2}{\sigma(q^k)} > \frac{8}{5}\cdot\left(\frac{n^2}{q^k}\right)$$
since $I(q^k) < 5/4$ holds unconditionally (i.e., for $k \geq 1$).  Additionally, note that
$$\frac{8}{5}\cdot\left(\frac{n^2}{q^k}\right) > \frac{8n}{5}$$
is true if $q^k < n$.  Furthermore, note that we then have the estimate $n > \sqrt[3]{N}$.

Lastly, note that we have
$$\frac{\sigma(n^2)}{q^k} = \frac{2n^2}{\sigma(q^k)} > 2n > \sigma(n)$$
if the stronger inequality $\sigma(q^k)<n$ holds.
\end{remark}

\section{The proof of Lemma \ref{lem:1}}
Let $N=q^k n^2$ be an odd perfect number with Euler prime $q$.  We want to show that the sum
$$\frac{\sigma(q^k)}{n}+\frac{\sigma(n)}{q^k}$$
is bounded from above if and only if the sum
$$\frac{q^k}{n}+\frac{n}{q^k}$$
is bounded from above.

To this end, note that we have the trivial inequalities
$$q^k < \sigma(q^k) < 2q^k$$
and
$$n < \sigma(n) < 2n$$
since both $q^k$ and $n$ are greater than one, and because $q^k$ and $n$ are deficient (being proper divisors of the perfect number $N=q^k n^2$).
These two sets of inequalities imply that
$$\frac{q^k}{n} < \frac{\sigma(q^k)}{n} < 2\cdot\frac{q^k}{n}$$
and
$$\frac{n}{q^k} < \frac{\sigma(n)}{q^k} < 2\cdot\frac{n}{q^k}$$
so that we obtain
$$\frac{q^k}{n}+\frac{n}{q^k} < \frac{\sigma(q^k)}{n}+\frac{\sigma(n)}{q^k} < 2\cdot\bigg(\frac{q^k}{n}+\frac{n}{q^k}\bigg).$$

First, we show that
$$\frac{\sigma(q^k)}{n}+\frac{\sigma(n)}{q^k} \text{ is bounded from above } \implies \frac{q^k}{n}+\frac{n}{q^k} \text{ is bounded from above }.$$
Suppose that
$$\frac{\sigma(q^k)}{n}+\frac{\sigma(n)}{q^k}$$
is bounded from above.  This implies that
$$\frac{\sigma(q^k)}{n}+\frac{\sigma(n)}{q^k} \leq C_1$$
for some absolute constant $C_1$.
But since
$$\frac{q^k}{n}+\frac{n}{q^k} < \frac{\sigma(q^k)}{n}+\frac{\sigma(n)}{q^k}$$
this implies that
$$\frac{q^k}{n}+\frac{n}{q^k} < \frac{\sigma(q^k)}{n}+\frac{\sigma(n)}{q^k} \leq C_1$$
which means that
$$\frac{q^k}{n}+\frac{n}{q^k} < C_1.$$
We conclude that
$$\frac{q^k}{n}+\frac{n}{q^k}$$
is bounded from above.

Next, we prove that
$$\frac{q^k}{n}+\frac{n}{q^k} \text{ is bounded from above } \implies \frac{\sigma(q^k)}{n}+\frac{\sigma(n)}{q^k} \text{ is bounded from above }.$$
Suppose that
$$\frac{q^k}{n}+\frac{n}{q^k}$$
is bounded from above.  This implies that
$$\frac{q^k}{n}+\frac{n}{q^k} \leq C_2$$
for some absolute constant $C_2$.
But since
$$\frac{\sigma(q^k)}{n}+\frac{\sigma(n)}{q^k} < 2\cdot\bigg(\frac{q^k}{n}+\frac{n}{q^k}\bigg)$$
this implies that
$$\frac{\sigma(q^k)}{n}+\frac{\sigma(n)}{q^k} < 2\cdot\bigg(\frac{q^k}{n}+\frac{n}{q^k}\bigg) \leq 2C_2$$
which means that
$$\frac{\sigma(q^k)}{n}+\frac{\sigma(n)}{q^k} < 2C_2.$$
We conclude that
$$\frac{\sigma(q^k)}{n}+\frac{\sigma(n)}{q^k}$$
is bounded from above.

This finishes the proof of Lemma \ref{lem:1}.

\begin{remark}
\label{rem:3}
In general, the function $f(z)=z+(1/z)$ is not bounded from above.  (To see why, it suffices to consider the cases $z \to 0^{+}$ and $z \to \infty$.)

This means that we do not expect the sum
$$\frac{\sigma(q^k)}{n}+\frac{\sigma(n)}{q^k}$$
to be bounded from above.
\end{remark}

\section{The proof of Theorem \ref{thm:1}}
Let $N=q^k n^2$ be an odd perfect number with Euler prime $q$.  We want to show that $\sigma(q^k) \neq \sigma(n)$.

Suppose to the contrary that $\sigma(q^k)=\sigma(n)$.  Then we obtain
$$\frac{\sigma(q^k)}{q^k}=\frac{\sigma(n)}{q^k}$$
and
$$\frac{\sigma(n)}{n}=\frac{\sigma(q^k)}{n}$$
from which it follows that
$$\frac{\sigma(q^k)}{n}+\frac{\sigma(n)}{q^k}=\frac{\sigma(q^k)}{q^k}+\frac{\sigma(n)}{n}=I(q^k)+I(n)<I(q^k)+I(n^2).$$
But Dris proved in \cite{Dris1} and \cite{Dris2} that
$$I(q^k)+I(n^2)<3$$
so that
$$\frac{\sigma(q^k)}{n}+\frac{\sigma(n)}{q^k}<3.$$
This means that
$$\frac{\sigma(q^k)}{n}+\frac{\sigma(n)}{q^k}$$
is bounded from above.  This contradicts Lemma \ref{lem:1} (see Remark \ref{rem:3}).

This finishes the proof of Theorem \ref{thm:1}.

\begin{remark}
\label{rem:4}
Similarly, we can show that $\sigma(n) \neq q^k$.  For suppose to the contrary that $\sigma(n)=q^k$.

Then we have
$$2>\frac{\sigma(q^k)}{n}\cdot\frac{\sigma(n)}{q^k}=\frac{\sigma(q^k)}{n}$$
since ${q^k}n$ is deficient (being a proper divisor of the perfect number $N=q^k n^2$).  But this implies that
$$\frac{\sigma(q^k)}{n}+\frac{\sigma(n)}{q^k}=\frac{\sigma(q^k)}{n}+1<3$$
from which it follows that
$$\frac{\sigma(q^k)}{n}+\frac{\sigma(n)}{q^k}$$
is bounded from above.  This contradicts Lemma \ref{lem:1} (see Remark \ref{rem:3}).
\end{remark}

\section{The proof of Lemma \ref{lem:3}}
Let $N=q^k n^2$ be an odd perfect number with Euler prime $q$.  We want to show that the inequality
$$\frac{\sigma(q^k)}{n}+\frac{\sigma(n)}{q^k}<I(q^k)+I(n)$$
holds if and only if the biconditional
$$q^k < n \iff  \sigma(n) < \sigma(q^k)$$
holds.

To this end, observe that we have the series of biconditionals
$$\frac{\sigma(q^k)}{n}+\frac{\sigma(n)}{q^k}<I(q^k)+I(n) \iff {q^k}\sigma(q^k)+n\sigma(n)<n\sigma(q^k)+{q^k}\sigma(n) \iff \bigg(q^k - n\bigg)\sigma(q^k)+\bigg(n - q^k\bigg)\sigma(n) < 0$$
$$\iff \bigg(q^k - n\bigg)\cdot\bigg(\sigma(q^k) - \sigma(n)\bigg) < 0 \iff \bigg(q^k < n \implies \sigma(n) < \sigma(q^k)\bigg) \land \bigg(n < q^k \implies \sigma(q^k) < \sigma(n)\bigg)$$ 
$$\iff \bigg(q^k < n \iff \sigma(n) < \sigma(q^k)\bigg).$$

Notice that we have used the facts that $q^k \neq n$ (since $\gcd(q,n)=1$) and $\sigma(q^k) \neq \sigma(n)$ (from Theorem \ref{thm:1}) as underlying assumptions throughout.

This finishes the proof of Lemma \ref{lem:3}.

\section{The proof of Corollary \ref{cor:1}}
Let $N=q^k n^2$ be an odd perfect number with Euler prime $q$.  We want to give an unconditional proof for the truth of the biconditional
$$q^k < n \iff \sigma(q^k) < \sigma(n) \iff \frac{\sigma(q^k)}{n} < \frac{\sigma(n)}{q^k}.$$
It suffices to show only the first biconditional
$$q^k < n \iff \sigma(q^k) < \sigma(n).$$
We consider three cases:

\textbf{Case 1}
$$\frac{\sigma(q^k)}{n}+\frac{\sigma(n)}{q^k}<I(q^k)+I(n)$$
We know (from \cite{Dris1} and \cite{Dris2}) that $I(q^k)+I(n)<I(q^k)+I(n^2)<3$.  This implies that
$$\frac{\sigma(q^k)}{n}+\frac{\sigma(n)}{q^k}$$
is bounded from above, which contradicts Lemma \ref{lem:1}.

\textbf{Case 2}
$$\frac{\sigma(q^k)}{n}+\frac{\sigma(n)}{q^k}=I(q^k)+I(n)$$
This equation is equivalent to
$$\bigg(q^k - n\bigg)\cdot\bigg(\sigma(q^k)-\sigma(n)\bigg)=0.$$
Since $q^k \neq n$, we must have $\sigma(q^k)=\sigma(n)$, contradicting Theorem \ref{thm:1}.

\textbf{Case 3}
$$\frac{\sigma(q^k)}{n}+\frac{\sigma(n)}{q^k}>I(q^k)+I(n)$$
This is equivalent to the inequality
$$\bigg(q^k - n\bigg)\cdot\bigg(\sigma(q^k)-\sigma(n)\bigg)>0,$$
which in turn is equivalent to the truth of the biconditional
$$q^k < n \iff \sigma(q^k) < \sigma(n).$$

This finishes the proof of Corollary \ref{cor:1}.

\section{Concluding Remarks}
Since ${q^k}n$ is deficient if $N=q^k n^2$ is an odd perfect number, then $I(q^k n)<2$.  This implies that
$$\frac{1}{2}\cdot\frac{\sigma(q^k)}{n}<\frac{q^k}{\sigma(n)}$$
and
$$\frac{1}{2}\cdot\frac{\sigma(n)}{q^k}<\frac{n}{\sigma(q^k)}$$
from which it follows that
$$\frac{1}{2}\cdot\bigg(\frac{\sigma(q^k)}{n}+\frac{\sigma(n)}{q^k}\bigg)<\frac{q^k}{\sigma(n)}+\frac{n}{\sigma(q^k)}.$$
Since the arithmetic mean is never less than the harmonic mean, and since
$$\frac{\sigma(q^k)}{n} \neq \frac{\sigma(n)}{q^k}$$
(see \cite{Dris3} for a proof of this inequation and some related considerations) then we have
$$\frac{2}{\frac{n}{\sigma(q^k)}+\frac{q^k}{\sigma(n)}}=\frac{2}{\frac{1}{\sigma(q^k)/n}+\frac{1}{\sigma(n)/q^k}} < \frac{1}{2}\cdot\bigg(\frac{\sigma(q^k)}{n}+\frac{\sigma(n)}{q^k}\bigg)$$
from which we obtain
$$\frac{2}{\frac{n}{\sigma(q^k)}+\frac{q^k}{\sigma(n)}}<\frac{q^k}{\sigma(n)}+\frac{n}{\sigma(q^k)}.$$
We conclude that
$$\sqrt{2}<\frac{q^k}{\sigma(n)}+\frac{n}{\sigma(q^k)}.$$
We now claim that either
$$\frac{\sigma(q^k)}{n}<\sqrt{2}<\frac{\sigma(n)}{q^k}$$
or
$$\frac{\sigma(n)}{q^k}<\sqrt{2}<\frac{\sigma(q^k)}{n}$$
holds.  (It suffices to prove one inequality, as the proof for the other one is very similar.)

To this end, assume that
$$\sqrt{2}<\frac{\sigma(n)}{q^k}.$$
This implies that
$$\sqrt{2}\cdot\frac{\sigma(q^k)}{n}<\frac{\sigma(q^k)}{n}\cdot\frac{\sigma(n)}{q^k}=I(q^k n)<2$$
which finally gives
$$\frac{\sigma(q^k)}{n}<\frac{2}{\sqrt{2}}=\sqrt{2}<\frac{\sigma(n)}{q^k}.$$
This proves our claim.

We now consider whether the following further refinements are possible:

\textbf{Case A}
$$1<\frac{\sigma(q^k)}{n}<\sqrt{2}<\frac{\sigma(n)}{q^k}<2$$
In this case,
$$\frac{\sigma(q^k)}{n}+\frac{\sigma(n)}{q^k}<2+\sqrt{2}$$
so that
$$\frac{\sigma(q^k)}{n}+\frac{\sigma(n)}{q^k}$$
is bounded from above.  This contradicts Lemma \ref{lem:1}.

\textbf{Case B}
$$1<\frac{\sigma(n)}{q^k}<\sqrt{2}<\frac{\sigma(q^k)}{n}<2$$
Similarly, in this case,
$$\frac{\sigma(q^k)}{n}+\frac{\sigma(n)}{q^k}<2+\sqrt{2}$$
so that
$$\frac{\sigma(q^k)}{n}+\frac{\sigma(n)}{q^k}$$
is bounded from above.  This contradicts Lemma \ref{lem:1}.

Consequently, since $\sigma(q^k) \neq n$ (because $\sigma(q^k) \equiv k+1 \equiv 2 \pmod 4$) and $\sigma(n) \neq q^k$, then we either have
$$\sigma(q^k)<n$$
or
$$\sigma(n)<q^k.$$

\begin{remark}
\label{rem:5}
The result in Corollary \ref{cor:1} together with the main findings in the preprint \cite{DagalDris} shows that
$$3 \nmid q^k n^2 \implies q < n.$$
This conclusion is derived independently of Brown's and Starni's methods.
\end{remark}

\begin{remark}
\label{rem:6}
By Corollary \ref{cor:1}, if $N=q^k n^2$ is an odd perfect number with Euler prime $q$, then there are a total of four cases to consider:
$$\textit{Case }\alpha: q^k < \sigma(q^k) < n < \sigma(n)$$
$$\textit{Case }\beta:  q^k < n < \sigma(q^k) < \sigma(n)$$
$$\textit{Case }\gamma: n < q^k < \sigma(n) < \sigma(q^k)$$
$$\textit{Case }\delta: n < \sigma(n) < q^k < \sigma(q^k)$$
Note that Cases $\beta$ and $\gamma$ imply that $k \neq 1$.  Also, from previous considerations, we know that $n < \sigma(q^k)$ and $q^k < \sigma(n)$ cannot both be true.  Consequently, Cases $\beta$ and $\gamma$ do not hold.

We are left with the scenarios:
$$\textit{Case }\alpha: q^k < \sigma(q^k) < n < \sigma(n)$$
$$\textit{Case }\delta: n < \sigma(n) < q^k < \sigma(q^k)$$

It turns out we can dispose of Case $\delta$ when $k=1$.  We obtain
$$\frac{\sigma(q^k)}{n}+\frac{\sigma(n)}{q^k}=\frac{\sigma(q)}{n}+\frac{\sigma(n)}{q}<\bigg(\sqrt{3}+(\sqrt[6]{3}\cdot{{10}^{-500}})\bigg)+1,$$
where the estimate
$$\frac{\sigma(q)}{n}<\sqrt{3}+(\sqrt[6]{3}\cdot{{10}^{-500}})$$
uses Acquaah and Konyagin's estimate $q<n\sqrt{3}$ \cite{AcquaahKonyagin} and Ochem and Rao's lower bound $N>{{10}^{1500}}$ for the magnitude of an odd perfect number \cite{OchemRao}. This implies that
$$\frac{\sigma(q^k)}{n}+\frac{\sigma(n)}{q^k}=\frac{\sigma(q)}{n}+\frac{\sigma(n)}{q}$$
is bounded from above, which contradicts Lemma \ref{lem:1}.

Consequently, $k \neq 1$ must hold in Case $\delta$.  From the papers \cite{Dris2} and \cite{Dris3}, this implies that $q < n$.

Since $\sigma(q^k) < n$ holds in Case $\alpha$, and since $q \leq q^k < \sigma(q^k)$, we also have $q < n$ under Case $\alpha$.
\end{remark}

We summarize the results we proved in Remark \ref{rem:6} in the following theorems.

\begin{theorem}
\label{thm:2}
If $N=q^k n^2$ is an odd perfect number with Euler prime $q$, then $q < n$ holds unconditionally.
\end{theorem}

\begin{theorem}
\label{thm:3}
If $N=q^k n^2$ is an odd perfect number with Euler prime $q$, then $k=1$ implies $\sigma(q^k)<n$.
\end{theorem}

\section{Further Research}
Let $N=q^k n^2$ be an odd perfect number with Euler prime $q$.  Suppose that the Descartes-Frenicle-Sorli conjecture that $k=1$ is true.

By Theorem \ref{thm:3} and Lemma \ref{lem:1}, $q+1=\sigma(q)=\sigma(q^k)<n$, so that we then have a further refinement of the following bounds (see the paper \cite{Dris3}):
$$\frac{\sigma(q)}{n}<1<I(q) \leq \frac{6}{5} < \bigg(\frac{5}{3}\bigg)^{\frac{\ln(4/3)}{\ln(13/9)}}<I(n)<2<\frac{\sigma(n)}{q}.$$

Again, by Lemma \ref{lem:1}, if $k=1$ then the ratio
$$\frac{\sigma(n)}{q^k}=\frac{\sigma(n)}{q}$$
is not bounded from above.  This implies that the ratio
$$\frac{\sigma(q^k)}{n}=\frac{\sigma(q)}{n}$$
is not bounded from below.  This means that we can take $\sigma(q)/n$ to be arbitrarily small, from which we conclude that $q$ has to be vastly smaller than $n$.

These considerations beg answers to several (obvious) questions, which we leave for other researchers to investigate.

\section{Acknowledgments}
The author thanks Carl Pomerance for sharing his expertise.  The author is also grateful to Keneth Adrian P. Dagal for helpful conversations that led to most of the results presented in this paper.  The author is also indebted to the anonymous referee(s) whose valuable feedback improved the overall presentation and style of this manuscript.

\makeatletter
\renewcommand{\@biblabel}[1]{[#1]\hfill}
\makeatother

\end{document}